\documentclass{mfatshort}
\begin{document}

\title[Topological structure of functions on a 3-manifold]
	{Topological structure of functions with isolated critical points on a 3-manifold}

\author{B. I. Hladysh}
\address{Department of Mechanics and Mathematics,
Taras Shevchenko National University of Kyiv,
Kyiv,Glushkova 4-e, 03127, UKRAINE}
%\curraddr{}
\email{bohdanahladysh@gmail.com, biv92@ukr.net}
\thanks{The second author was partially supported by the Austrian Academy of Sciences in the frame of the Project between the Austrian Academy of Sciences and the National Academy of Sciences of Ukraine on New mathematical methods in astroparticle and quantum physics.}

\author{A. O. Prishlyak}
\address{Department of Mechanics and Mathematics,
Taras Shevchenko National University of Kyiv,
Kyiv,Glushkova 4-e, 03127, UKRAINE}
\email{prishlyak@yahoo.com}
%\thanks{Thank you very much.}

%    General info
\subjclass[2010]{57R45, 57R70, 58C27} % e.g. 35A30; 81Q05
% For 2010 Mathematics Subject Classification see http://www.ams.org/mathscinet/msc/msc2010.html
%\date{DD/MM/2004}
%\dedicatory{This paper is dedicated to you.}
\keywords{Topological equivalence, critical point, 3-manifold}

\begin{abstract}
To each isolated critical point of a smooth function on a 3-manifold we put in correspondence a tree (graph without cycles). We will prove that functions are topologically equivalent in the neighborhoods of critical points if and only if the corresponding trees are isomorphic. A complete topological invariant of functions with isolated critical points, on a closed 3-manifold, will be constructed.
\end{abstract}

\maketitle

\section{Introduction} 		      % don't type final punctuation

There are many works focused on topological properties of functions defined on manifolds. The first ones in such area were Kronrod's \cite{Kronrod} and Reeb's \cite{Reeb} papers. Let $M$ be a smooth 3-manifold and $f, g: M~\to~\mathbb{R}$ be smooth functions.

The functions $f$ and $g$ are called \emph{topologically equivalent}, if there are homeomorphisms $h : M \to M $ and $k : \mathbb{R} \to \mathbb{R}$ such that $f\circ h = k\circ g$.

We say that functions are \emph{topologically conjugated}, if they are topologically equivalent and homeomorphism $k$ preserves orientation. In this case $h$ and $k$ will be called by \emph{conjugated homeomorphisms}.

Notice that functions without critical points can be topologically equivalent to function with critical points. For example, the function $f(x,y)=x^2+y$ doesn't have critical points, but the function $g(x,y)=x^2+y^3$ has a critical point --- $(0,0)$ and topological equivalence of these functions can be defined in the following way: $h(x,y)=(x,y^3), k=id_{\mathbb{R}}.$

A critical point of a function $f$ is \emph{removable} if the function $f$ is topologically equivalent to the function without critical points in some neighborhoods of this point.

Differentiable equivalence is studied in the theory of singularities. It is such a topological equivalence in which the conjugated homeomorphisms are diffeomorphisms.

The problem of topological classification of Morse functions was solved in \cite{Arn}, \cite[p.49--129]{BF}, [5 -- 9], \cite[p.19--48]{Sha90} for closed manifolds of different dimensions. The same result for arbitrary functions with isolated critical point on closed 2-manifolds was obtained in \cite{Pri02}. The relevance of this problem is contributed by the close connection with the Hamiltonian dynamical system's classification in dimensions 2 and 4.

In this paper we give a local topological classification of functions with isolated critical points and a global topological classification of smooth functions with three critical points on closed 3-manifolds.

Takens \cite{Tak68} has proven that an isolated critical point of a smooth function on 3-manifold have conic type. We use this result for local topological classification of functions. We also construct a colored graph that classifies functions with 3 critical points on 3-manifolds. We note that each 3-manifold admitting a function with 3 critical points is a connected sum of several copies of $S^1 \times S^2$ in oriented case or $S^1 \widetilde{\times} S^2$ (the non-trivial fiber bundle over $S^1$ with the fiber $S^2$) in non-oriented case.

\section{Topological structure of a neighborhood of a critical point}

Let $f$ be a smooth function on a smooth 3-manifold $M$. It is known \cite{Pri02} that if $p$ is an isolated critical point and  $y=f(p)$, then there exists a closed neighborhood $U (p)$ such that $$f^{- 1 } (y )  \cap   U (p) = Con ( \cup S _ {i}^{1 } ), $$ where $ Con (\cup S_{i}^{1 } )$ is a cone over a disjoint union of circles $S_{i}^{1 } $, that is the union of two-dimensional disks, whose centers are pasted together into the point $p$.

In order to describe the behavior of a function in a neighborhood of a critical point $p$ we  will construct a tree (graph without cycles) $Gf_{p}$. Let $U(p)$ be the neighborhood described above, whoce boundary is a sphere $S^{{ 2}} $ and $\partial  (f^{- 1 } (y )  \cap   U  ( p )) = \cup S_{i}^{1 }$ is the union of the embedded circles. To each component $D_{j} $ of $S^{{ 2}} \setminus \cup S_{i}^{1 } $ we put in correspondence a vertex $v_{j} $ of the graph  $Gf_{p} $ and to each circle $S_{i}^{1 } $ we associate an edge $e_{i} $. The vertex $v_{j}$ is incident to $e_{i}$ if and only if the boundary of $D_{j} $ contains  $S_i^1$. Thus, $v_i$ and $v_j$ are connected by an edge if $D_i$ and $D_{j}$ are neighbor.

\textbf{Example 2.1.}
The function $ f (x, y, z) = x ^ {2} + y ^ {2} -z ^ {2} $ has the number of circles $k =2$, and the function $ f (x, y, z) = (x ^ {2} + y ^ {2} -z ^ {2}) (x ^ {2} + y ^ {2} -4z ^ {2} ) (x ^ {2} -4y ^ {2} + z ^ {2}) $ has the number of circles $k = 6$ at the critical point $p$, that is the origin. The location of these circles on the sphere, as well as the corresponding graphs $Gf_p$ is shown in Fig. 1. (The sphere is regarded as a plane with a point at infinity).

\begin{figure}[h]
\centering
\includegraphics[width=3.99in,height=1.60in]{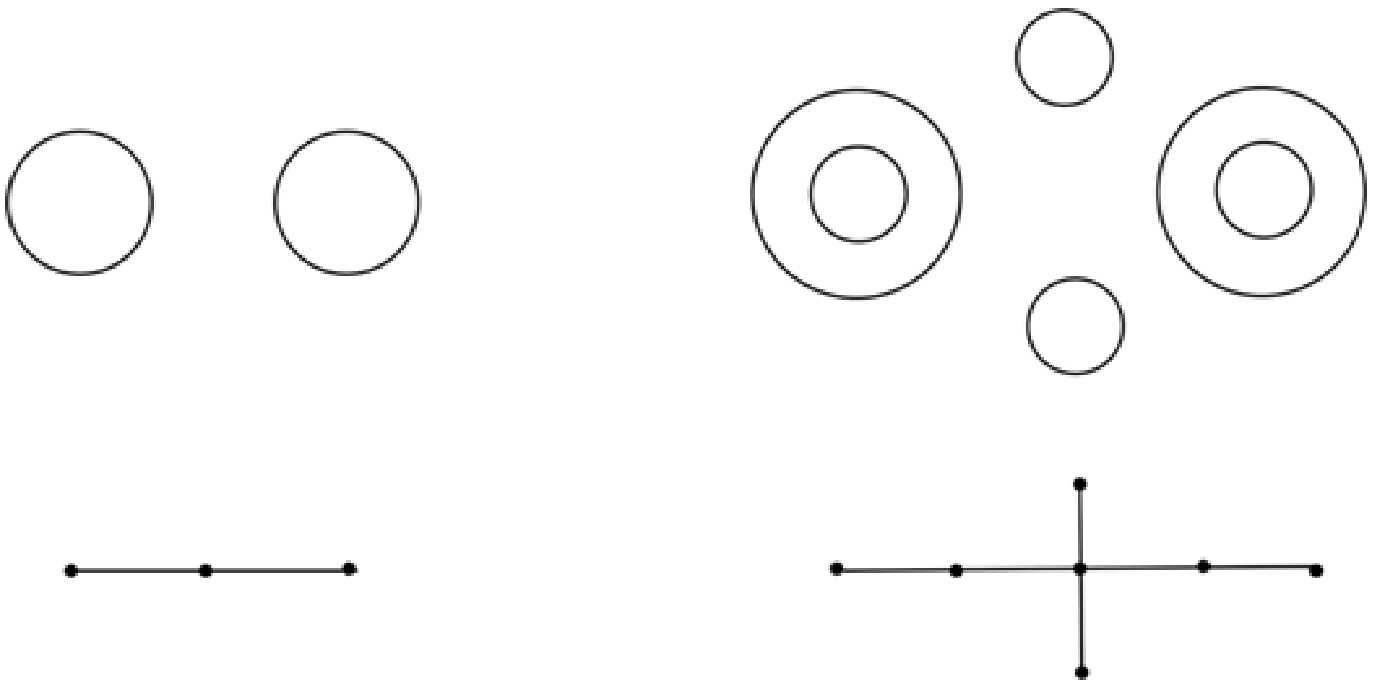}

Fig. 1.
\end{figure}

For a function $f$ and its isolated critical point $p$ we define a neighborhood $W_{p}$ of a point $p$ in $f^{-1}(f(p))$, which is homeomorphic to $Con ( \cup S_{i}^{1 } )$. Following Takens [11] let $$ W_{p} ( \varepsilon ) = \{ x \in M : \vert f  ( x ) - y _0 \vert < \varepsilon , cl ( \gamma (x )) \cap W_{p} \ne \emptyset \}$$ be a neighborhood $W_{p}(\varepsilon )$ of $p$ in $\mathbb{R}^{{ 3}} $ for $\varepsilon
> 0$, where $ \gamma  (x )$ is the integrated trajectory of the gradient field of $f$ which contains the point $x$. The above neighborhood $W_{p}(\varepsilon )$  will be called \emph{canonical}.

Let $ W_{p} ( \varepsilon )$ be a canonical neighborhood then it is the boundary $ \partial W_{p} $ ($ \varepsilon $) is a union: $$ \partial   W_{p} ( \varepsilon ) = V_- \cup V_{{ 0}} \cup V_{+} ,$$  where $V_- = \partial W_{p}  ( \varepsilon ) \cap f^{- 1 }  (f  ( p ) - \varepsilon ),$  $V _ {+} = \partial W_{p}  ( \varepsilon )  \cap f^{- 1 }  (f  ( p ) + \varepsilon ),$ $V_{-} = cl ( \partial W_{p}  ( \varepsilon )  \setminus  (V_ - \cup V_{+} )) = cl( \partial W_{p} ( \varepsilon )  \cap f^{- 1 } (f  ( p ) - \varepsilon $, $f  ( p ) + \varepsilon $)). We will say that $V_{+} $ and $V_- $ are the upper and the lower foundation, accordingly, and $V_{{ 0}} $  is the side wall of the neighborhood $W_{p} $ ($ \varepsilon $).

The side wall is a union of closed neighborhoods of circles $S_{i}^{1 } $. Therefore $ V_{{ 0}}  =  \cup S _{i}^1 \times  [- \varepsilon , \varepsilon ].$

Thus for every $i $ and $t \in  [- \varepsilon, \varepsilon]$: $S_{i}^{1 } \times \{ t \} \subset f^{-1}(f(p)+t),$ and for each $s \in S_{i}^{1 } $: $ s\times  [-\varepsilon ,  \varepsilon ]  \subset \gamma  (s, 0).$

On each cylinder $S_{i}^{1 } \times  [- \varepsilon , \varepsilon]$ level lines of the function and the integral curves define a structure of a direct product.

We denote by $D_{i}' $ a connected component of $V_{+}$ or $V_{-}$ which belongs to $D_{i}$. Then $$D_{i}'  = D_{i} \setminus \cup S_{i}^{1 } \times (- \varepsilon , \varepsilon ).$$

Let $W_{i} = \mbox{Con}(S_{i}^{1} )  \subset f^{- 1 }  (f  ( p ))$ and $ D_{i}'' $ be the set of those points from $D_{i} ' $ whose the integrated trajectories have $p$ as a limit point: $$D_{i} '' = \{x \in D' \vert   \gamma  (x ) \cap W = \emptyset\}.$$ Then $D_{i} '' $ is a deformation retract of both $D_{i} ' $ and $D_{i} $ (see example 2.3).

For the determinacy let $f  ( D_{i} '' ) = f  ( p ) + \varepsilon$. We will construct a new vector field $X'$ on the set $$ U_{i} '  = \{x \in cl (W_{p}  ( \varepsilon )) \vert \gamma  (x ) \cap W \ne \emptyset , f  ( x ) \geq f ( p ) \}  \cong  (0,1] \times S_{i}^{1} \times  [0, \varepsilon ].$$  In order to do this we will consider coordinates $(u, s, t )$ on $U_{i} ' $. Coordinate $t$ of the points $x$ is equal to $f(x)-f(p)$. Since $Con(S_{i}^{1}) \setminus p $ is homeomorphic to $(0,1] \times S_{i}^{1 }$, coordinates $u $ and $s$ at $t = 0$ are defined by that homeomorphism. For an arbitrary point $x \in U_{i} ' $ we will choose coordinates $u$ and $s$ to be equal to coordinates $u $ and $s $ of the point $ \gamma (x)\cap cl (W_{p} )$. The existence of such coordinates follows from the tube theorem for flows (or about a rectification of a vector field). Since the integral curves of a vector field grad$f$ coincide with coordinate lines, it has coordinates $\{0, 0, v(x)\}$, where $v(x)>0$ for each point $x \in U_{i} $.

Let $$X' = \{ \frac {u \cdot v (x)}{\sqrt {u^{2} + t^{2}}} , \frac {t \cdot v (x)} {\sqrt {u^{2} + t^{2}} }, 0 \} \ \ \ \ \ \mbox{ if } t \geq 2 u\varepsilon;$$
$$X'=\{ {\frac {{(1 - u) \cdot v (x)}} {{\sqrt {(1 - u)^{2} + (2\varepsilon - t)^{2}}}}} {, \frac {{(2\varepsilon - t) \cdot v (x)}} {{\sqrt {(1 - u)^{2} + (2\varepsilon - t)^{2}}}}}, 0 \} \ \ \ \ \ \mbox{ if }   t \leq       2 \varepsilon u \Large{\}}.$$

It is easy to see that $X' $ is a gradient-like field for $f $, and $X'$ coincides with $X$ at points with coordinate $u = 1$ and at $u \to  0$. The set $D_{i} ' \setminus D_{i} '' $ for $X'$ consists of points with coordinate $u  > 1/2$. The vector field $X'$ that is constructed in a such way  is called an \emph{inclined} vector field. In contrast to gradient fields which depends on Riemannian metric, all inclined vector field of a function $f$ are topologically equivalent.

\textbf{Theorem 2.1.}
Let $p$ and $q$ be isolated critical points of smooth functions $f: \mathbb{R}^{3} \to \mathbb{R}^{1}$ and $g: \mathbb{R}^{3} \to \mathbb{R}^{1}$ correspondingly. Then there are neighborhoods $U$ of $p$ and $V$  of $q$ and homeomorphisms $h: U \to V$ and $k : \mathbb{R} \to \mathbb{R}$ such that $f\circ h = k\circ g$ if and only if graphs $Gf_{p} $ and $Gg_{q} $ are isomorphic.

\emph{Proof.}
\textit{Necessity.} It follows from the construction of the graphs that the restriction of a homeomorphism $h $ on the boundaries of the neighborhoods will determine a required isomorphism of the graphs.

\textit{Sufficiency.} Fix an isomorphism $i: Gf_{p} \to Gg_{q}$. Let $W_{p}(\varepsilon )$ be a canonic neighborhood of point $p$ and $ \pi : W_{p}  ( \varepsilon )  \to W_{p} $ be a map given by the formula $$ \pi (x) = \left\{
\begin{array}{cc}
p,  & \mbox{ if } \ \ \gamma (x) \cap f^{- 1} (f (p)) = \emptyset,
\hfill
\\  \gamma (x) \cap f^{- 1} (f (p)),    \ &
\mbox{ if } \ \
 \gamma (x) \cap f^{- 1} (f (p)) \ne \emptyset. \hfill \\
\end {array} \right. $$

For $q$ and its canonical neighborhood $W_{q} $ ($ \varepsilon $) define $ \pi $ in a similar way. Let us construct a homeomorphism of boundary sphere $H $: $\partial W_{p} $ ($ \varepsilon $) $ \to   \partial W_{q} $ ($ \varepsilon $) such that at each point $x \in W_{p} $ ($ \varepsilon $): $$\vert f(x)-f(p) \vert  =  \vert g(H(x))- g( q )\vert.$$

We construct a required homeomorphism of boundary spheres. The isomorphism of the graph sets correspondence of $D_{j}$ of two functions and also correspondence of $S_{i}^{1}$. We choose arbitrary orientation of the edges of one of the graphs and orient the edges of another graph in a  such way that isomorphism of the graphs preserves the orientation. Fix also orientations of the spheres. Then the orientation of the edges of the graphs determines the orientation of the circles. Fix an arbitrary homeomorphisms  $\varphi_{i} $ of the  circles to the relevant circles that preserve the orientations. These homeomorphisms multiplied by identical map of a segment $[- \varepsilon ,  \varepsilon ]$ define homeomorphisms of cylinders. The products $ \varphi_{{ i}} $ with identical map of $(0,1]$ define homeomorphisms of cones Con($S_{i}^{1} $). We can extend homeomorphisms of cones on sets $U_{i} '$ using inclined vector field and coordinates that are relevant to them. It follows from the construction of inclined fields that they define homeomorphisms of boundaries of region $U_{i} $ (on equality of the relevant coordinates). Thus we have homeomorphisms of boundaries of $D_{i} '' $. We extend them inside of $D_{i} '' $ arbitrarily. These homeomorphisms define correspondences of integrated trajectories. The correspondence of points of trajectories is given by equality of a difference of values of functions in them with a value at critical points.

It follows from the construction that the constructed map is a homeomorphism which maps levels of the function $f$ into levels of the function $g$.

\textbf{Example 2.2.}
Accordingly to the notations, described in the proof of Theorem 2.2 for the function $ f(x,y,z) = x^2 - y^2 + z^2$ in standard Riemannian metric, the neighborhood $W_p (\epsilon )$  is shown in Fig. 2. Wherein $V_- = {D}_1 '  \cup D_3'$, $V_+ = {D}_2'$, ${D}_1'' ={p}_ 1$, ${D}_2 '' = {S}^1$, ${D}_3 {''} = {p}_ 3$.

\begin{figure}[h]
\centering
\includegraphics[width=3.99in,height=2.00in]{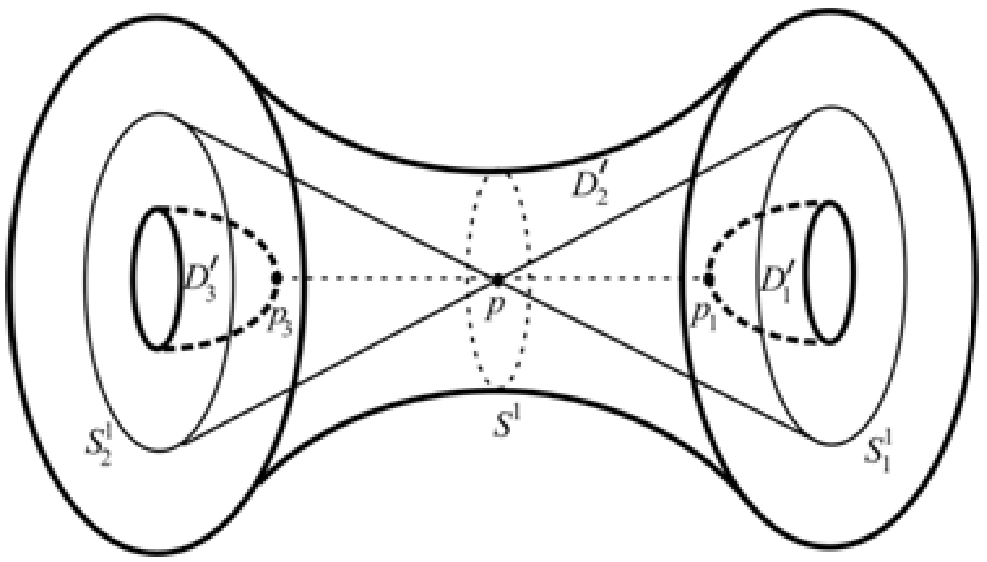}

Fig. 2.
\end{figure}

\section{Global conjugation of functions with tree critical points on closed 3-manifolds}

In order to solve the problem of global conjugation of functions with three critical points on 3-manifolds we will construct a decomposition of manifold into 3-dimensional disks intersected on boundaries. Thus each disk will contain a unique critical point and is isotopic to the neighborhood, which was constructed in section 2. The global invariant of a function consists of these disks and maps with a glued boundaries.

Further we will consider functions on a closed 3-manifold. Such a function has at least two critical points --- minima and maxima. If a function defined on a closed manifold has exactly two critical points, then according to Reeb theorem that manifold is a sphere and critical points are points of minima and maxima. This situation is not interesting.

%\subsection{Distinguishing graph of a function %with three critical points.}

Let $M$ be a closed 3-manifold and $f : M \to  \mathbb{R}$ be a smooth function with three isolated critical points $p_1, p_2, p_3$ and critical values $y_{i}=f(p_i), i=1,2,3$ such that $y_1 < y_2 <y_3$.

Denote by $U_{2} = W_{p_{2}} ( \varepsilon )$ the neighborhood of the point $p_2$ constructed in section 1. Let $U_{1 }$ be the connected  component of $M\setminus (cl (U_2)  \cup f^{- 1 } (p_2))$ which contains the point $p_{1 } $ and $U_3$ is one for the point $p_3 $. Then each $U_{i}$ is homeomorphic to the open 3-disk. Indeed, $U_2$ is homeomorphic to 3-dimentional disc according to the construction of $W_{p_{2}} ( \varepsilon )$ in the previous section, $U_1$ includes single critical point (minima). That is why $U_1$ can be considered as 0-handle $h^0=D^0\times D^3,$ which is homeomorphic to $D^3.$ In the same way $U_3$ includes a single critical point (maxima) and it can be considered as $h^3=D^3\times D^0.$ Thus, $U_3$ is homeomorphic to $D^3$ \cite[p.75]{Matsumoto}. We consider boundary spheres $S_i^2  = \partial U_{i}, i=1,2,3 $. We color components of $S_1^2 \cap S_2^2$ in the white color, components of $S_{{ 2}}^{{ 2}} \cap S _ {{ 3}}^{{ 2}} $ in the black color and  components of $S_{1 }^{{ 2}} \cap S_{{ 3}}^{{ 2}} $ in the gray color. Thus each of spheres $S_{i}^{2} $ is decomposed into parts of two colors. As well as in section 2 on each sphere we construct the graphs  $Gf_{i} $. We color vertices of the graphs in colors of the corresponding parts of $S_{i}^{2}$. At each edge of each graph fix a new vertex and split that edge into two subedges. The obtained graphs will be denoted by $Gf_{i}'$. From these graphs we paste together the new graphs $Gf$ in such a manner that two old vertices of different graphs are glued together if the region relevant to them coincide in $M $. Also we  paste together incident edges if the circles that are relevant to them coincide. In another way  graphs $Gf$ can be obtained if we regard component of $S_{1 }^{{ 2}} \cap S_{{ 2}}^{{ 2}} \cap S_{{ 3}}^{{ 2}} $ and component of   $S_{1 }^{{ 2}} \cup S_{{ 2}}^{{ 2}} \cup S_{{ 3}}^{{ 2}} \setminus S_{1 }^{{ 2}} \cap S_{{ 2}}^{{ 2}} \cap S_{{ 3}}^{{ 2}} $ as vertices. Thus two vertices are connected by an edge if one of them corresponds to component $U_i$ of $S_{1 }^{{ 2}} \cup S_{{ 2}}^{{ 2}} \cup S_{{ 3}}^{{ 2}} \setminus S_{1 }^{{ 2}} \cap S_{{ 2}}^{{ 2}} \cap S_{{ 3}}^{{ 2}} $ and another one to the component of it boundary $\partial U_i$. We color vertices relevant to the 2-dimensional components in the same colors as the components. Fix the orientation on these spheres. If we get the same orientation on each component of spheres intersection, we say that the orientation number of the vertex is $-1$, and the orientation number is $+1$ if the orientations are different. In such a way each white vertex gets orientation number ($+1$ or $-1$). The graph $Gf$ has following properties:

1) The vertices of the graph are divided into four types: white, black, gray and non-colored. The number of vertices of each color (the first three types) is same. The non-colored vertices have degree 3. Each white vertex is equipped with the orientation number ( $+1$ or $-1$).

2) If from the graph we remove vertices of one color and edges that incident to them, we obtain simply-connected graphs (tree) $Gf_{i}'$.

\textbf{Definition 3.1} A graph will be called by \emph{distinguishing} if it satisfies two properties described above. Two distinguishing graphs will be called \emph{equivalent} if there is an isomorphism of one them to another, which maps vertices of each color onto vertices of the same color and preserve the orientation numbers %$+1$ and $-1$
in each white vertex or changes all these numbers.  The distinguishing graph constructed as above for the given function will be called by \emph{distinguishing graph of} this \emph{function}.

\textbf{Example 3.1}  We will construct a function $f$ with three critical points  and the diagram for $f$ on the $S^1 \times S^2$. Heegaard diagram of this manifold is a torus,  each meridian system of which consists of a closed curve and the two curves do not intersect (see. Fig. 3).

\begin{figure}[h]
\centering
\includegraphics[width=2.5in,height=1.3in]{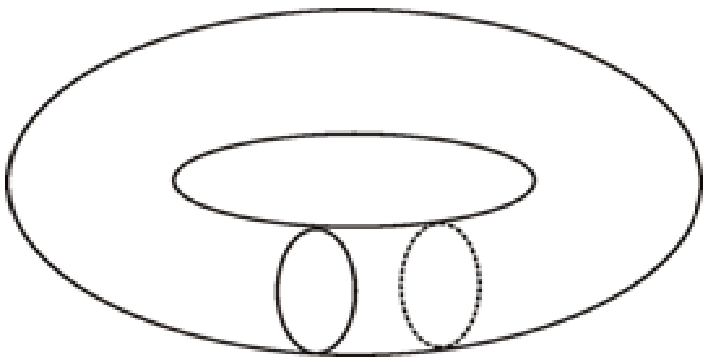}

Fig. 3.
\end{figure}

Since the meridians do not intersect each other, there is a handle decomposition which has one handle of each index, and the handle of index 1 and 2 are attaching at the same time. Thus in Fig. 4 it is shown a 2-handle $h^2=D^2\times D^1$ in the form of thickened hemisphere. This hemisphere is glued according to the embedding from the $S^1\times D^1$ to the boundary of a 3-body. Also it is described 1-handle $h^1=D^1\times D^2$ which is glued according to the embedding of $S^0\times D^2$ into a 3-body, such that one of components lays inside the image of the $S^1\times D^1$ and another one is outside this image.

\begin{figure}[h]
\centering
\includegraphics[width=3.0in,height=1.50in]{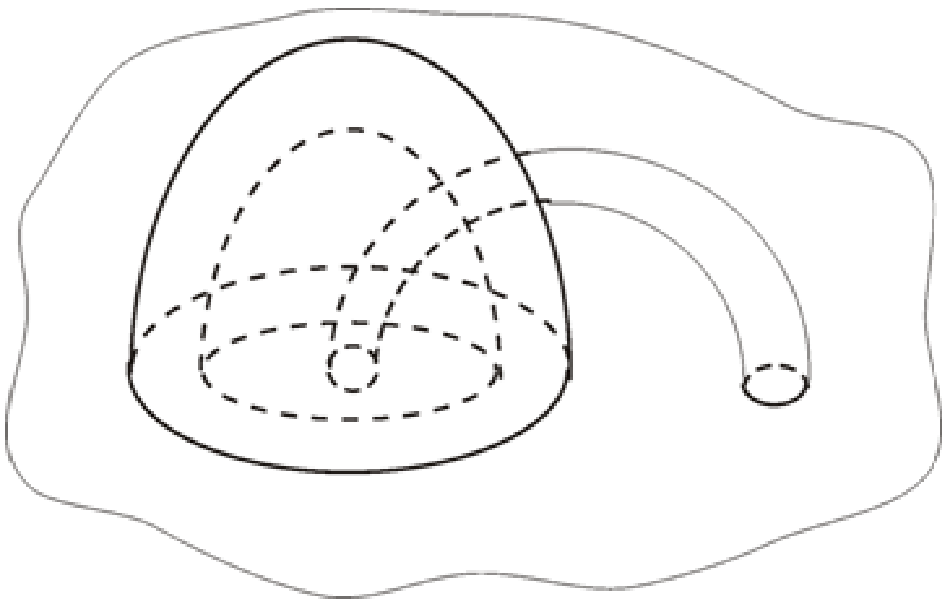}

Fig. 4.
\end{figure}

Gluing of the two handles may be replaced by gluing of a 3-dimensional disc as in Fig. 5. In this figure one of the ends of the 1-handle is glued inside the image of $S^1\times D^1$ and another one is glued outside 2-handle.

\begin{figure}[h]
\centering
\includegraphics[width=3.0in,height=1.50in]{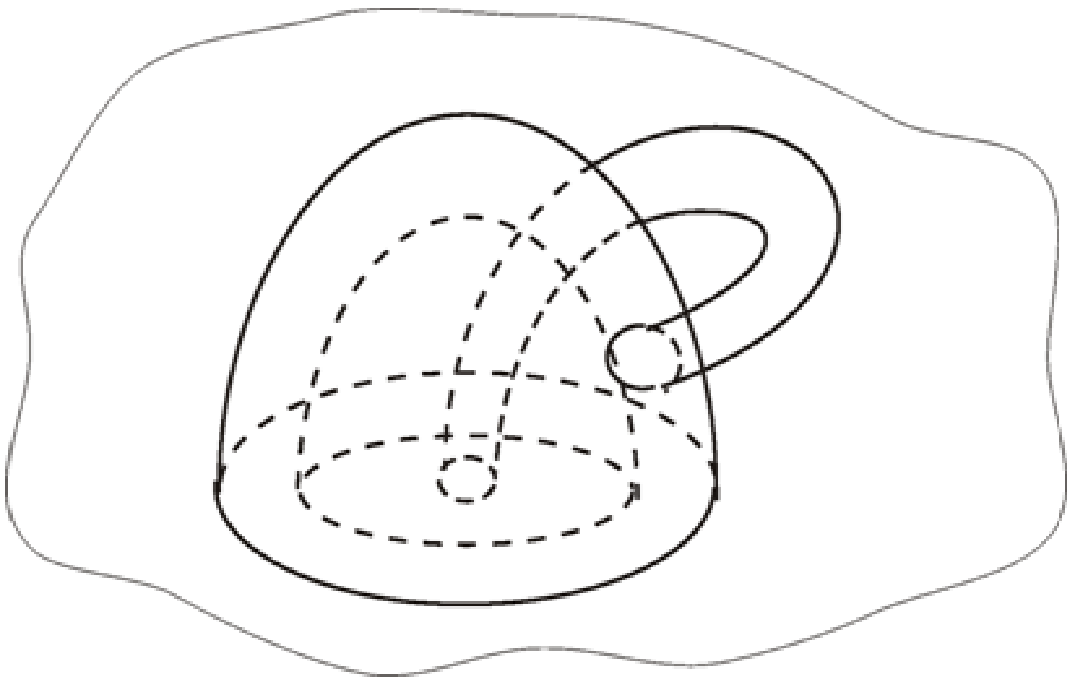}

Fig. 5.
\end{figure}

Thus, we have a decomposition of  $S^1 \times S^2$ into 3 three-dimensional disk: 0-handle, 1-handle $ \cup $ 2-handle (as in Fig. 7) 3-handle. The second disc is attached to the first by an embedding whose image is shown in Fig. 6.

\begin{figure}[h]
\centering
\includegraphics[width=2.0in,height=1.00in]{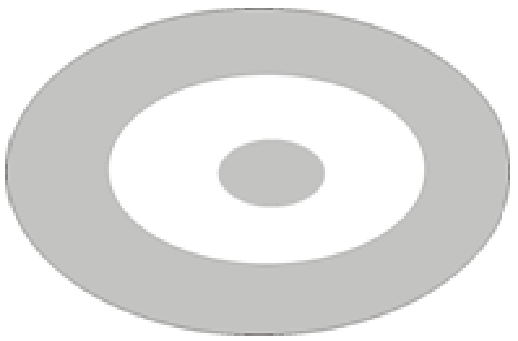}

Fig. 6.
\end{figure}

Regions on the boundaries of other discs (the boundaries of the second and third discs) have the same form. Then all the graphs of critical points (matching color) have the form shown in Fig. 7. Both white vertices have the orientation numbers $+1$.

\begin{figure}[h]
\centering
\includegraphics[width=2.5in,height=0.16in]{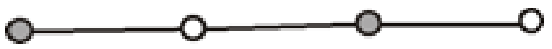}

Fig. 7.
\end{figure}

After splitting each edge in half and gluing respective halves of different graphs we obtain the graph shown in Fig. 8.

\begin{figure}[h]
\centering
\includegraphics[width=1.7in]{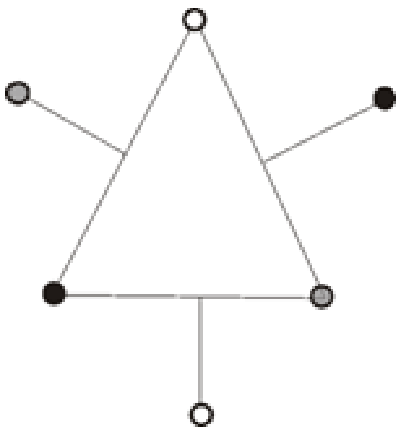}

Fig. 8.
\end{figure}

This graph is the graph of functions with three critical points on $S^1 \times S^2$.

\textbf{Lemma 3.1}
The orientation numbers ($+1$ or $-1$) of white vertices define the signs of another vertices up to the multiplier $\pm 1$ of all these vertices.

\emph{Proof.}
Suppose that the orientation numbers of white vertices ($+1$ or $-1$) are given. Consider any black vertex and set its orientation number to $+1$. Further  define the orientation numbers of the nearest (corresponding to this black vertex) gray vertices. Let us consider one of them. There is exactly one $\mathrm{T}$-vertex (uncolored vertex with valency 3) between above-described gray vertex and the black one. Note that the vertices being adjacent with the $\mathrm{T}$-vertex either all have orientation numbers $+1$ or two of them have the orientation number $-1$ and one has the orientation number $+1$. It is the result of the following considerations: if we define the orientation at the point of intersection of three 3-dimension disks (spheres), then we get three orientation numbers $+1$ in the cases of their coincidence if they have opposite orientations, and the changing of the orientation of one of these disks leads to the changing of two signs (two intersection regions). Thus, only the following combinations of orientation numbers are possible: $+1, -1, -1$; $+1, -1, -1$ and $+1, +1, +1$.

The orientation number of a gray vertex is determined by the signs of black and white vertices being adjacent to the previously described $\mathrm{T}$-vertex (as a product of these numbers). In particular, if a black vertex has the orientation number $+1$, then the gray  and the white vertices will have the same number. Then, by the same arguments, the signs of black vertices being adjacent to the described above gray one are determined from the sign of this gray vertex (if these signs are not defined earlier). Note that the gray-black subgraph is a tree. That is why there is a single shortest path from fixed black vertex to any other vertex and as a result the orientation number of a black vertex defines the numbers of the others gray and black vertices. Also changing the number of black vertex leads to the changing of the numbers of each black and gray vertices simultaneously. Lemma 1 is proved.

\textbf{Theorem 3.1.}
Let  $f, g: M \to \mathbb{R}$ be smooth functions which have three critical points on a smooth closed 3-manifold $M$. The functions $f$ and $g$ are conjugated if and only if their distinguishing graphs are equivalent.

\emph{Proof.}
\textit{Necessity.} The conjugated homeomorphism sets a homeomorphism of $U_{i} $, $i = 1,2,3$  and  induces an isomorphism between distinguishing graphs.

\textit{Sufficiency.} Not losing generality we can assume that the functions have same critical values $-1$, 0 and 1. Suppose that distinguishing graphs of functions are isomorphic. After removal of gray vertices and incident to them edges and replacement of non-colored vertices with two incident edges by one edge we will receive isomorphic graphs. This allows us to construct a homeomorphism $h $ of region $U_2$. The correspondence between gray vertices defines a correspondence between connected components of  $f^{- 1 }  (0)  \setminus U_{{ 2}} $ and $g^{-1} (0) \setminus h  ( U_{ 2} )$. All these components are oriented surfaces of genus 0 (because they are subsets of the 2-sphere). It follows from isomorphism of the distinguishing graphs that corresponding $U_i$ have the same number of boundary components (equal to the degree of the gray vertex). Therefore, the homeomorphisms of the edges, given by $h $, can be extended to a homeomorphism $h_{{ 0}} $ of $f^{- 1 } $ (0) $ \setminus U_{{2}} $. Consider gradient vector fields $\mbox{grad} (f)$ and $\mbox{grad} (g)$ in the Riemannian metric used in the construction of $h $. The homeomorphisms $h $ and $h_{{ 0}} $ induce a bijective correspondence of trajectories of the gradient vector fields. The homeomorphisms of the relevant trajectories, which preserve the values of the function, gives a required homeomorphism of the functions.

\textbf{Notation 3.1.}
Two functions are conjugate if their graphs are equivalent or will be equivalent in the result of colors replacement of white vertexes to black and black to white.

Denote by $n$ the number of the  vertexes colored in one of three color on the distinguishing graph. In this case we will say that $n$ is the \emph{complexity} of the corresponding function.

\textbf{Notation 3.2.} In oriented case $M$ is homeomorphic to the connected sum  $\sharp_n S^1\times S^2$ (where $n$ is the complexity number) and in non-oriented case $M$ is homeomorphic to $\sharp_n S^1\widetilde{\times}S^2$.

If the manifold is oriented, then there are a unique distinguishing graph if $n = 1$ and a unique distinguishing graph if $n = 2$ (up to topological equivalence). Thus, there is a unique (up to topological equivalence or conjugation) function of complexity 1 (on $S^3$), and a unique function of complexity 2 (on $S^1\times S^2$). Their distinguishing graphs are shown in Fig. 9. In this case the last graph has the similar orientation number of white vertices.  If the orientation numbers of white vertices are different (non-oriented manifold), we get a unique (up to topological equivalence) function defined on $S^1\widetilde{\times}S^2$.

\begin{figure}[h]
\centering
\includegraphics[width=3.6in]{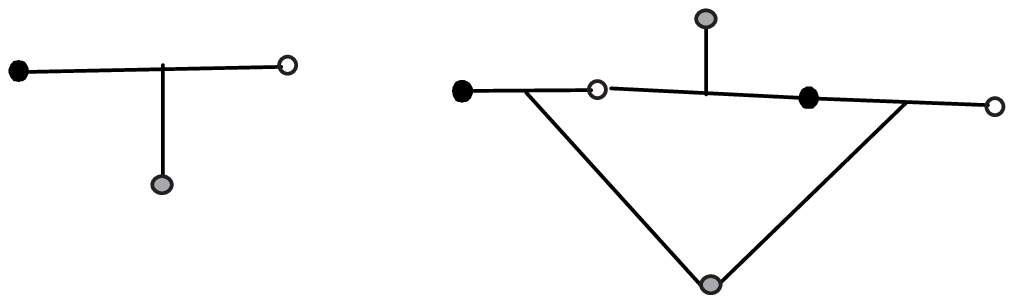}

Fig. 9.
\end{figure}

%There are six  topologically not equivalent functions if $n  = 3$. Their distinguishing graph are figured in the fig. 10. Thus if in the graphs 4), 5) and 6) one exchange color of white vertexes to black and black to white, then the resulting distinguishing graphs are not equivalent to original ones. For the graphs 1), 2) and 3) such replacement of the graphs leads to graph that equivalent original ones. Thus, there exists 9 topologically not conjugate functions of the complexity 3.

%\begin {figure}[h]
%\centering
% \includegraphics[width=4.0in]{s10.eps}

%Fig. 9.
%\end {figure}

\begin {figure}[h]
\centering
 \includegraphics[width=3.8in]{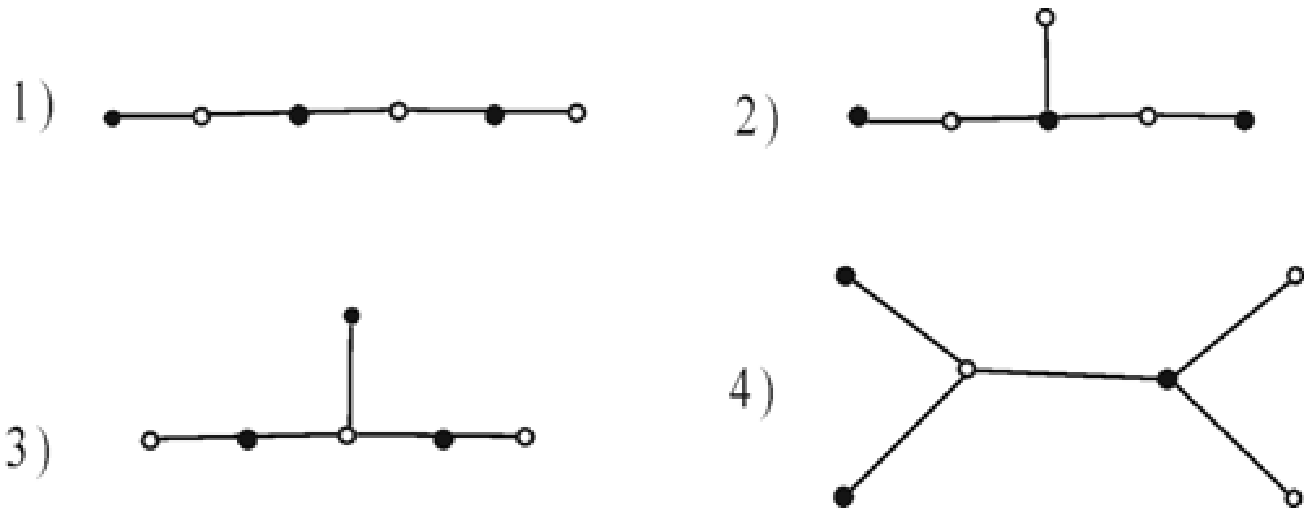}

Fig. 10.
\end {figure}

\begin {figure}[h]
\centering
 \includegraphics[width=4.0in,height=3.50in]{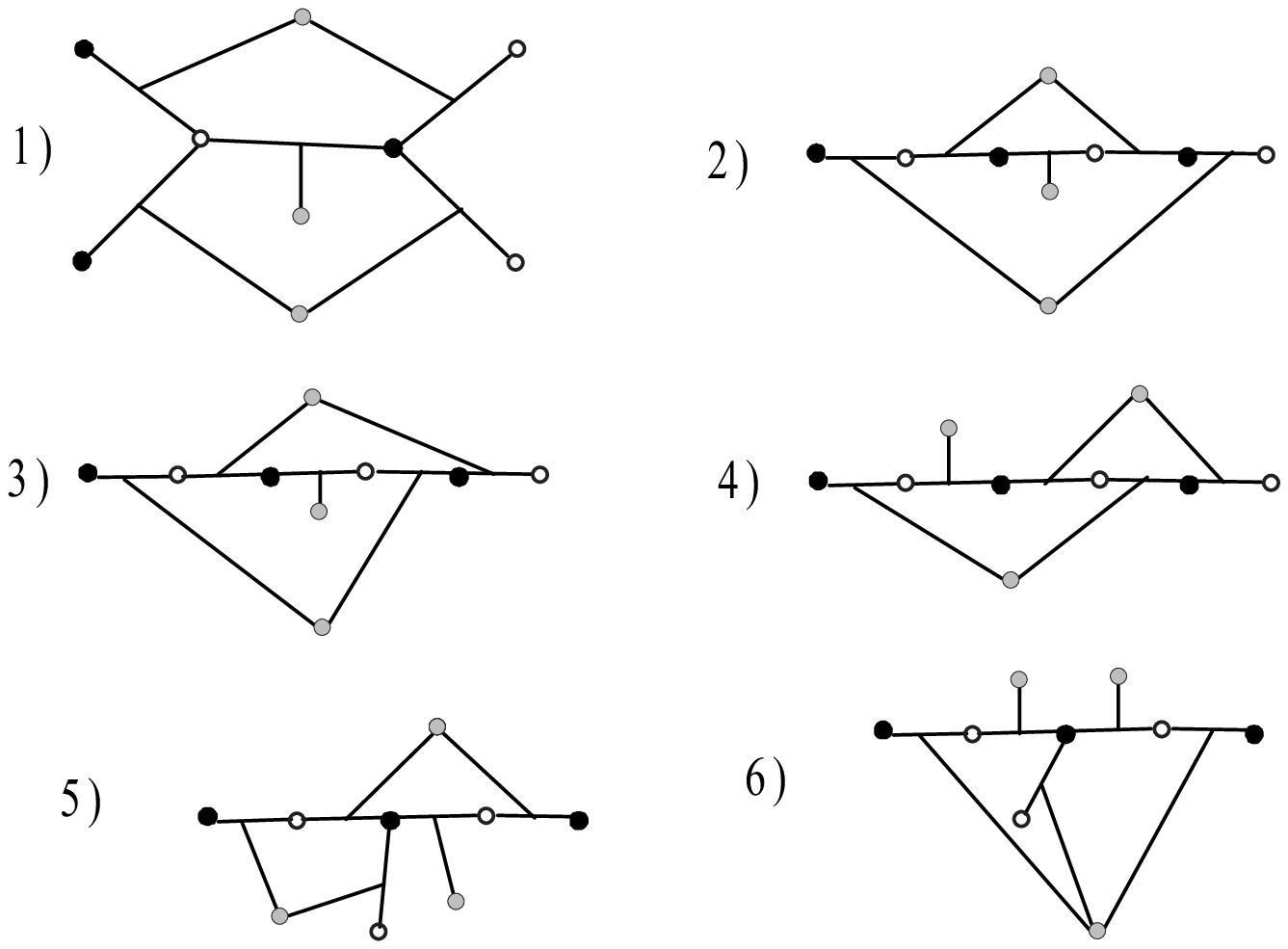}

Fig. 11.
\end {figure}

%\begin{center}

%\begin {figure}[h]
%\includegraphics[width=3.8in]{s12.eps}
%\includegraphics[width=3.8in]{s13.eps}

%\centering
%Fig. 11.
%\end {figure}
%\end{center}

%\begin {figure}[h]
%\centering
% \includegraphics[width=4.0in]{s13.eps}
%
%Fig. 12.
%\end {figure}

%\begin {figure}[h]
%\centering
% \includegraphics[width=4.0in,height=3.50in]{s14.eps}

%Fig. 13.
%\end {figure}

\begin{center}

\begin {figure}[h]
\includegraphics[width=3.9in]{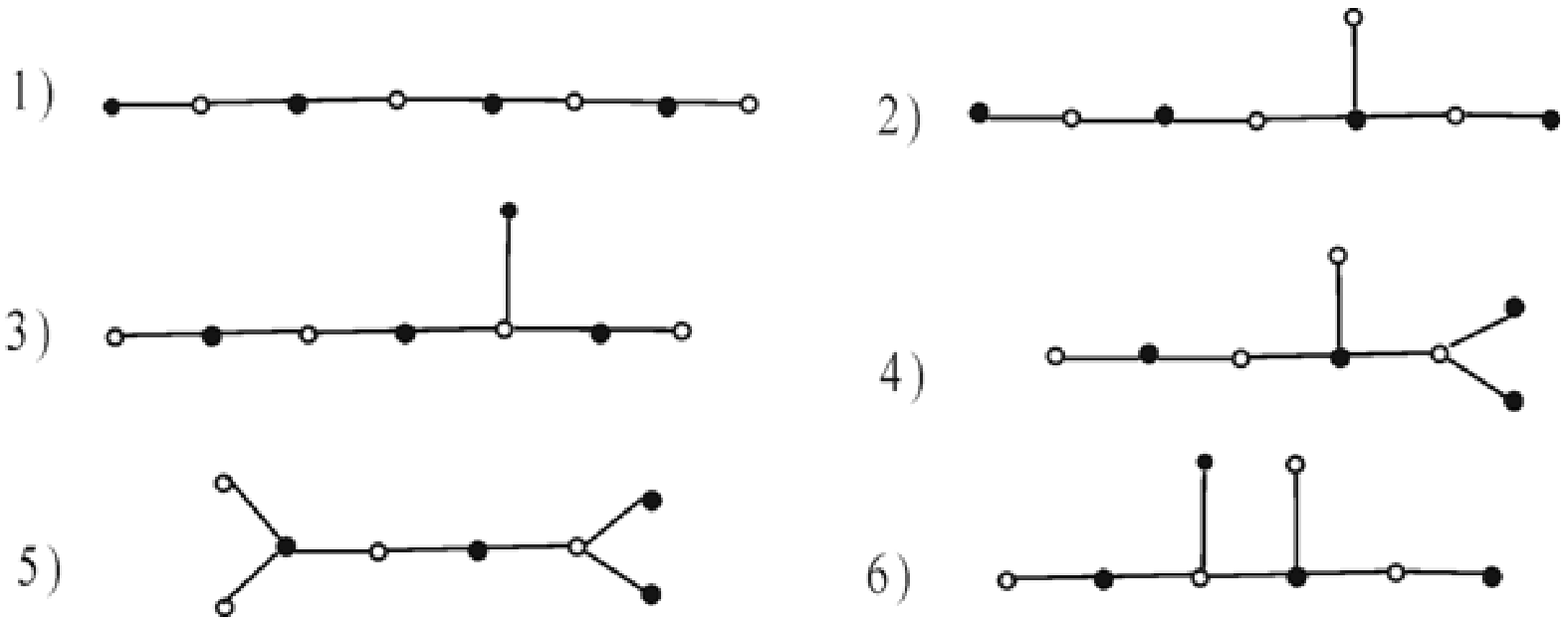}
\includegraphics[width=3.9in]{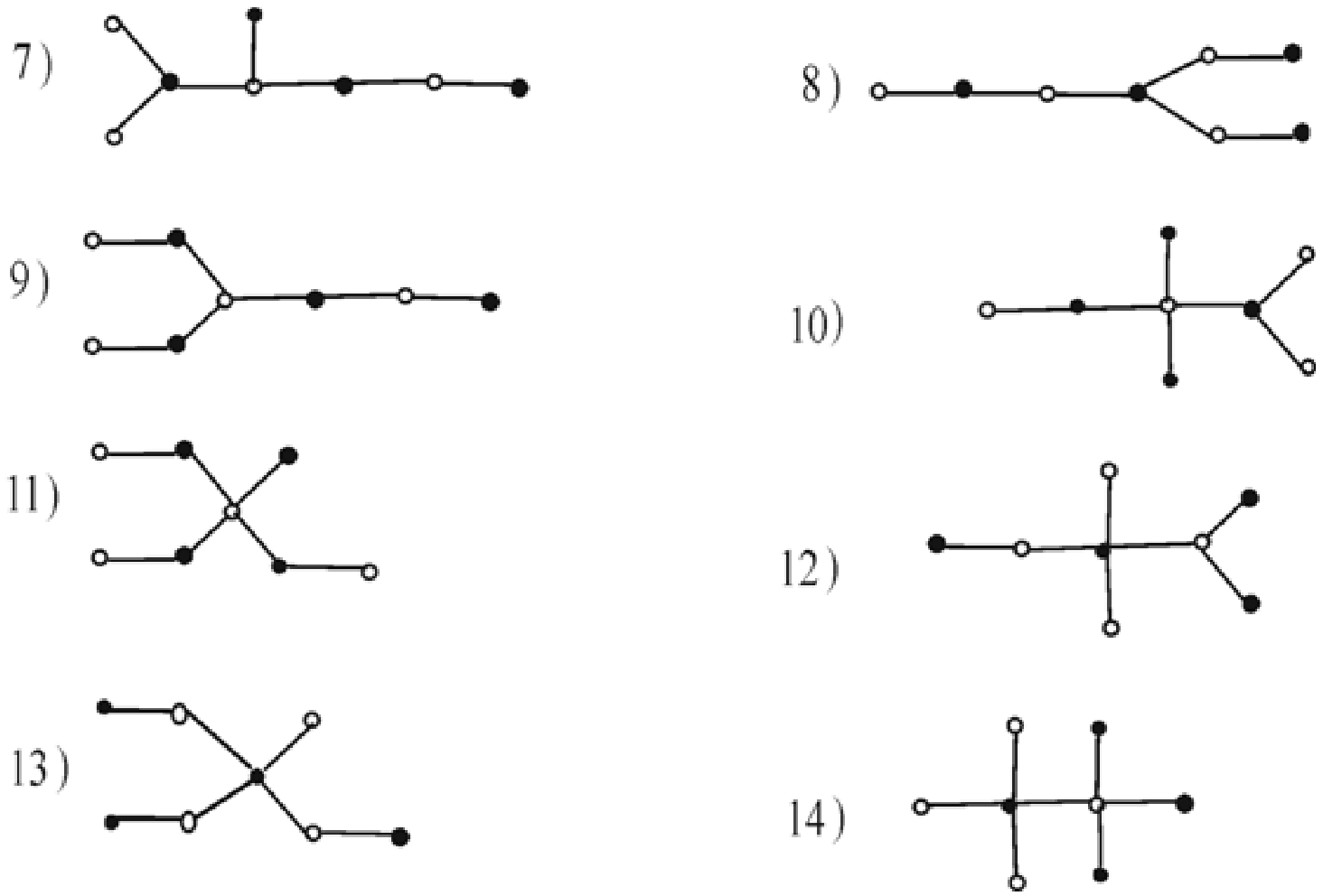}

\centering
Fig. 12.
\end {figure}
\end{center}

There are four (based on coloring) graphs of critical points if $n  = 3$.  They are shown in Fig. 10. The graph of function is obtained from the two graphs of points (one with white-black vertices, and other with white-gray vertices) by gluing homeomorphism of the neighborhood of the corresponding white vertices.  Obviously, we can glue  graphs having the same valence of corresponding white vertices.

There are six topologically non-equivalent functions if $n = 3$ defined on $(S^1\times S^2)\sharp (S^1\times S^2)$. Their distinguishing graph are shown in Fig. 11. Thus if in the graphs 4), 5) and 6) one exchange color of white vertexes to black and black to white, then the resulting distinguishing graphs are not equivalent to original ones. For the graphs 1), 2) and 3) such replacement of the graphs leads to graph which is equivalent to original ones. Thus, there exists 9 topologically non-conjugate functions of complexity 3.

In the case of non-oriented manifold $(S^1\widetilde{\times}S^2)\sharp(S^1\widetilde{\times}S^2)$ there exist two non-equivalent sign sets of white vertices, see graphs 1) and 6), Fig. 11. In the same way we get three sign sets of white vertices of graphs 2) -- 5) of Fig. 11. Thus, there are 16 non-equivalent functions with three critical points on $(S^1\widetilde{\times}S^2)\sharp(S^1\widetilde{\times}S^2)$. Concerning the conjugation of functions, we get twice as many non-conjugated functions corresponding to graphs 4) -- 6), Fig. 11, because the change of the color in each of these graphs leads to two different graphs. In such a way, we get  24 non-conjugated functions on $(S^1\widetilde{\times}S^2)\sharp(S^1\widetilde{\times}S^2)$.

In the case of complexity $n=4$ $(S^1\times S^2)\sharp(S^1\times S^2)\sharp(S^1\times S^2)$, we can get 14 (subject to coloring) graphs of critical points. They are shown in Fig. 12.

The numbers of admissible gluing of graph pairs are shown in Table. 1.

The sum of all the numbers in Table is equal to 179. So there is 179 topologically non-conjugate functions of complexity 4.

\begin{table}

\begin{tabular}{|p{0.15in}|p{0.15in}|p{0.15in}|p{0.15in}|p{0.15in}|p{0.15in}|p{0.15in}|p{0.15in}|p{0.15in}|p{0.15in}|p{0.15in}|p{0.15in}|p{0.15in}|p{0.15in}|p{0.15in}|} \hline
\textbf {} & \textbf {1} & \textbf {2} & \textbf {3} & \textbf {4} & \textbf {5} & \textbf {6} & \textbf {7} & \textbf { 8} & \textbf {9} & \textbf {10} & \textbf {11} & \textbf {12} & \textbf {13} & \textbf {14} \\\hline
\textbf{1} & 20 & 12 & & & & & &7 & &  & & & 1 & \\\hline
\textbf{2}  & 12 &8 & & & & & &5 & & & &  & 1 &  \\\hline
\textbf{3} & & & 2 & 1 & 1 & 8 & 3 &  &6 & & & 1 & & \\\hline
\textbf{4}  & & &1 & & & 1 & 2 &  &3 & & & & & \\\hline
\textbf{5} & & &1 & & & 1 &  &  &3 & & & & & \\\hline
\textbf{6} & & & 8 & 1 & 1 & 4 & 2 &  & 1 & & &1 & & \\\hline
\textbf{7} & & & 3& 2 &  & 2 & & & & & & & & \\\hline
\textbf{8}  & 7 &5 & & & & & & 4& & & &  & 1 &  \\\hline
\textbf{9} & & & 6 & 3 & 3 & 1  & & & 1 & & & 1 & & \\\hline
\textbf{10} & & & & & & & & & & & 2 & & & \\\hline
\textbf{11} & & & & & & & & & & 2 & 1  & & &1 \\\hline
\textbf{12} & & &1 & & &1 & & &1 &   & & & &  \\\hline
\textbf{13} &1 & 1 & & & & & &1 & & & &  & 1 &  \\\hline
\textbf{14} & & & & & & & & & & & 1 & & & \\\hline
\end {tabular}

\caption{\textbf {Number of different allowable gluing pairs of graphs}}
\end{table}

Gluing of graphs of type 1 can be made with permutations. Thus, at each vertex of valence 3 we have two numbers: 1) its sequence number in black and white graphs; 2) its sequence number in gray and white graphs if you start counting from the black (gray) vertex of the valence 1. Then there is a permutations of 7 numbers. But as the last number is always equal to 7, then it is given by the permutations of 6 numbers. On the other hand, this permutations indicates how the neighborhoods of white vertices glue, and hence defines the graph function. There are 20 permutations that define permissible graphs:

%\begin {enumerate}
   (5, 6, 3, 4, 2, 1),
   (4, 3, 6, 5, 2, 1),
   (3, 4, 6, 5, 2, 1),
   (4, 3, 5, 6, 2, 1),

   (5, 6, 4, 3, 1, 2),
   (2, 1, 5, 6, 4, 3),
   (6, 5, 2, 1, 4, 3),
   (6, 5, 1, 2, 4, 3),

   (5, 6, 1, 2, 4, 3),
   (6, 5, 2, 1, 3, 4),
   (5, 6, 2, 1, 3, 4),
   (2, 1, 6, 5, 3, 4),

   (2, 1, 5, 6, 3, 4),
   (2, 1, 4, 3, 6, 5),
   (3, 4, 2, 1, 6, 5),
   (4, 3, 1, 2, 6, 5),

   (3, 4, 1, 2, 6, 5),
   (6, 5, 3, 4, 1, 2),
   (3, 4, 6, 5, 1, 2),
   (4, 3, 5, 6, 1, 2).
%\end {enumerate}

If we allow substitution of black and gray vertices between them, then the number on the diagonal can be changed, and the numbers of gluing graphs of two different types should be considered once (instead of two, as was done above). There are twelve different 3-graphs obtained by gluing two graphs of type 1 (equivalent pairs of permutations are 1--18, 2--7, 3--8, 4--10, 6--12, 9--19, 11--20, 15--16, substituting 5, 13, 14, 17 will go into themselves), six 3-graphs of type 2, one 3-graphs of type 3, three 3-graphs of type 6 and eight 3-graphs of type 3. To sum up all the numbers, we find that the number of topologically non-equivalent functions of complexity 4 defined on $(S^1\times S^2)\sharp(S^1\times S^2)\sharp(S^1\times S^2)$ equals 93.

\section{Functions with 4 critical points}

Let $M$ be a closed oriented 3-manifold and $f : M \to \Re $ be a smooth function with 4 isolated critical points $p_1, p_2, p_3, p_4 $ and correspondent critical values $y_{i}=f(x_i), i=1,2,3,4 $ such that $y_{i} < y_{j}$, if $i < j. $ Thus $p_{1 } $ is a minimum point and $p_{{ 4}} $ is a maximum point.

We denote by $U_{i} = W_{p_{}} { (} \varepsilon {)} $ a neighborhoods of the points $p_{{ i}}$, i=2,3, defined in section 2. Let $U_{1} $ be a connected component of  $M\setminus $ (cl ($U_{{ 2}} $) $ \cup f ^ {- 1 } $ ($p_{{ 2}} $)), containing a point $p_{1}$, and  $U_{4} $ be a component of $M\setminus $ (cl ($U_{{ 3}} $) $ \cup f ^ {- 1 } $ ($p_{{ 3}} $)), containing a point $p_{4}$. Consider surfaces $F =\partial $ ($U_{1 } \cup U_{{ 2}} $) and $F ' = \partial $ ($U_{{ 3}} \cup U_{{ 4}} $). According to the construction they are homeomorphic to a regular level $f^{ - 1 } (z )$, where $y_2 < z < y_3$. Moreover, part $M_{{ 0}} $ of the manifold $M $, which located between them is homeomorphic to the cylinder $F\times $ [0,1]. Denote by \{$u_{i} $ \} closed curves of $cl(U_{1 })$ $ \cap $cl($U_{{ 2}} $) $ \cap M_{{ 0}} $ and by \{$v_{i} ' $ \} of $cl(U_3) \cap cl (U_4) \cap M_0$. Let $\pi$ be the projection of top base of cylinder $M_0$ to bottom base and let $v _ {i} = \pi $ ($v_{i} ' $). The surface $F$ constructed in such manner together with two sets of closed curves \{$u_{i} $ \}, \{$v_{i} $ \} on it will be called by \emph{a diagram} of the function $f$ and designated by $D_{f}$. While it is possible to describe gluing $U_{{ 2}} $ to $U_{1 } $ and $U_{{ 3}} $ to $U_{{ 4}}$ using graphs $G_{1 } $ and $G_{{ 2}} $, we will need the diagram of function to give the attaching $U_{{ 2}} \cup U_{1 } $ to $U _ {{ 3}} \cup U_{{ 4}} $.

The diagram of a function is similar to the Heegaard diagram of 3-dimensional manifold and for it one can use concepts of isomorphism, isotopy and semiisotopy of the diagrams \cite{Pri98}.

In the construction of the diagram of a function we have an ambiguity in the choice of a structure of a direct product on $M_{{ 0}} $. Thus the change of the structure of a direct product induces an isotopy of curves $v_{i} $, leaving curves $u_{i} $ invariant. On the contrary, each isotopy induces change from fixed structure to a new structure of a direct product. Thus the obtained diagrams are semiisotopic. Using semiisotopy (choices of a structure of a direct product on $M_{{ 0}} $) we cancel all twoangles in the diagram and  we  obtain the normalized diagram.

Two curves $u_{i} $ and $v_{j} $ are called parallel if they are isotopic in the complement to other curves, i.e. if they form the boundary of a connected component homeomorphic to $S^{1 } \times $ [0,1], obtained by splitting $F $ by curves $u_{i} $ and $v_{j} $.

\textbf{Proposition 4.1.} Two normalized diagrams are semiisotopic if and only if one can pass from the first diagram to the second one by isotopies consisting of permutations of parallel curves.
%iff from first diagram it is possible to obtain diagram, isotopic second one using permutations of parallel curves.

\emph{Proof.}
If diagrams do not contain parallel curves, the proof coincides with the proof of the similar statement for Heegaard diagrams \cite{Pri98}. If the parallel curves intersect in two points, then they form two twoangles. Depending on the way of reducing the twoangle we can obtain two distinct normalized diagrams. These diagrams differ by permutation of two parallel curves. The inverse is obvious: the permutation of two parallel curves can be obtained using the semiisotopy.

Consider the graphs $G_{1 } $ defining gluing of $U_{{ 2}} $ to $U_{1 } $. Each non colored vertex of valence 3 ( we denote by $V_{{ 0}} $ the set of such vertices) corresponds to a meridian from the first system of meridians. Vertices colored in the first and third colors (we denote the sets of such vertices $V_{1 } $ and $V_{{ 3}} $, accordingly), correspond to components into which the first system of meridians divides a surface. Thus we set a bijective map $ \psi_{1}: \{V_0,V_1 \cup V_2  \}  \to  \{\{u_{i}  \},  \pi_0 (F\setminus\cup u_{i} ) \}$. Similar correspondences $ \psi_{{ 2}}$ arise for the graph $G_{{ 2}} $ and the second system of meridians.

\emph{A scheme} of function $f $ is a quintuple $\{D_{f}, G_{1}, \psi_{1 },  G_{2},  \psi_{2} \}$, consisting of the diagram $D_{f}$ of the function, two graphs $G_1, G_2$ and maps of correspondence $\psi_{1 }, \psi_{2}$. Two schemes are called  \emph{equivalent} if there are isomorphisms of the diagrams and the graphs, compounded with maps of correspondence.

\textbf{Theorem 4.1.} Functions $f$ and $g$ with 4 critical points on 3-manifolds $M, N$ correspondingly are topologically conjugate if and only if the scheme of one of them is equivalent to the scheme obtained from another one by a semiisotopy of its diagram.

\emph{Proof.}
\textit{Necessity} follows from the construction and previous arguments.

\textit{Sufficiency.} Without loss of generality we may assume that the functions have critical values - 2, - 1, 2, 3. As well as in Theorem 3.4 we construct a homeomorphism from $f^{-1 } ([-2,0])$ to $g^{-1}([-2,0])$ and a homeomorphism from $f^{-1}([1,3])$ to $f^{-1} ([1,3])$. Thus on sets $f^{-1}([0,1])$ and $g^{-1} ([0,1])$ the structures of a direct product are fixed which were used in construction of the diagrams. Then the semiisotopy of the diagrams can be considered as a level-by-level homeomorphism from $f^{-1}([0,1])$ to $g^{{ -} 1 } $ ([0,1]). On requirements of the theorem it coincides on the boundaries with the constructed above homeomorphisms from $f^{-1} ([-2,0])$ and $f ^ { -1} ([1,3])$ and thus it is extension of these homeomorphisms up to a required homeomorphism of manifold.

\section{Function with an arbitrary finite number of critical points}

Let $p_{1 } $, \ldots, $p_{k} $ be critical points, $f (p_{1}) \leq f  (p_2) \leq  \ldots  \leq f ( p_{k} )$. We fix a Riemannian metric on the manifold and small enough closed neighborhoods $W_{1 } $, \ldots, $W_{k} $ \textit  of critical points (not intersected among themselves and having the same structure as $W_{p} $ ($ \varepsilon $) in Theorem 2.2). We want to construct analogue of a handle decomposition that is the decomposition of the manifold into a union $M=H_{1 } \cup $ \ldots $ \cup H_{k} $, where $W_{i} \subset H_{i}, i = $ 1, \ldots, $k $. We construct neighborhoods $H_{1 } $, \ldots, $H_{k} $ by an induction and call them generalized handles. Put $H_{1 } $ = $W_{1 } $. Let $S$ ($ W_{i} $) be a set of those points, whose positive orbit with respect to the gradient field grad $f $ intersects $W_{i} $. Then we put $H_{i} $ = cl ($S $ ($ W_{i} $) $ \setminus \cup_{j { <} i} H_{j} $). Thus, $M $ can be obtained by a sequential gluing of generalized handles.

Consider the following set $S =\partial H_{1 } \cup $ \ldots $ \cup \partial H_{k} $. It has natural structure of the stratified set. Thus each strata of dimensionality 2 lays in the intersection of two different generalized handles $H_{j} $, $H_{k}$. By the \emph{diagram} of a function we will mean the stratified set, constructed on it, together with a pair of the numbers of adjacent strata for every strata and sets of the numbers having the same critical values. As well as before the diagram will be called \emph{homeomorphic}, if there is a homeomorphism of the stratified sets preserving the pairs and the sets of the numbers.

\textbf{Proposition 5.1.} Two functions with isolated critical points on 3-manifolds are topologically conjugate if and only if it is possible to construct their diagram that is homeomorphic.

\emph{Proof.}
\textit{Necessity.} The restriction of a conjugated homeomorphism on the first stratified set induces a homeomorphism between the first diagram and the diagram constructed on the image this map.

\textit{Sufficiency.} As well as above without losing generality it is possible to assume that the functions have the same sets of critical values.

The boundary of each generalized handle can be divided into three parts: 1) lower base consisting of boundary intersections with handles, which have the smaller numbers, 2) upper base consisting of points, in which gradient field is transversal to boundary of the handle and not included in the lower base, 3) side walls consisting of points not inherings to the bases. The given homeomorphism of the stratified sets can easily be improves so that it will amp side walls on side walls, and preserve their partition on levels of the function. Then Theorem 2.2 allows us to extend this homeomorphism up to a required conjugated homeomorphism.

Let us consider a problem when not homeomorphic diagram correspond to topologically conjugate function. Similarly to functions with 4 critical points a choice of another Riemannian metric implies in that generalized handles have isotopic attaching maps. The correspondent diagrams will be called \textit{semiisotopic}.

\textbf{Proposition 5.2.}
Two functions with isolated critical points on 3-manifolds are topologically conjugate if and only if their diagrams are semiisotopic.

The proof is similar to the proof of Theorem 4.2.

\textbf{Remark 5.1.}
The obtained criteria of topological conjugations are not constructive in the cases of 3 and 4 critical points .

\end{document}